\documentclass[12pt]{article}
\usepackage{epsfig}
\usepackage{amsmath}
\usepackage{amsthm, amssymb, amscd}

\author{Hugh Howards,
John Luecke
}

\title{Strongly n-trivial Knots}
\date{}

\begin{document}
\maketitle

\newtheorem{definition}{Definition}[section]
\newtheorem{corollary}[definition]{Corollary}
\newtheorem{theorem}[definition]{Theorem}
\newtheorem{lemma}[definition]{Lemma}
\newtheorem{claim}[definition]{Claim}
\newtheorem{ex}[definition]{Example}
\newtheorem{question}[definition]{Question}
\newtheorem{exer}[definition]{Exercise}
\newtheorem{conjecture}[definition]{Conjecture}

\begin{abstract}
We prove that given any knot $k$ of  genus $g$,
 $k$ fails to be strongly $n$-trivial for all $n$,
$n \geq 3g-1 $. 

\medskip


\medskip


\medskip
 
\end{abstract}


\section{Introduction}

We start with a little background.

\begin{definition}
A knot $k$ is called  ``$(n$-$1)$-trivial''
 if there exists a projection of $k$,
such that one can choose $n$ disjoint sets of
crossings of the projection with the property 
that  making the
crossing changes corresponding to any of the $2^{n}-1$
nontrivial combination of the sets of crossings
turns the original knot into the
unknot. The collection of sets of 
crossing changes is said to be an ``$(n$-$1)$-trivializer for $k$''.
\end{definition}

\begin{conjecture}
The unknot is the only knot that is $n$-trivial for all $n$. 
\label{conj:setsarcs}
\end{conjecture}
{\em Note:} \hspace{.5mm} A knot that is $n$-trivial is 
automatically $m$-trivial for all $m \leq n$.

\medskip
 Work of Gusarov
[Gu] and Ng and Stanford [NS] shows that this question equates to showing
that the only knot with vanishing Vassiliev invariants for all $n$ is the 
unknot.  Thus, Conjecture~\ref{conj:setsarcs} is at the heart of the study of
Vassiliev invariants.

One reason why this question is interesting is that 
it takes a geometric
approach to Vassiliev invariants, instead of the traditional 
algebraic approach and therefore is relatively unexplored.  
Vassiliev invariants measure geometric 
properties of knots, which in turn are geometric objects,
so it is reasonable to hope that the geometry 
might play an integral role
in their study.  

The following definition derives 
its motivation from $n$-triviality.

\begin{definition} 
A knot $k$ is called ``strongly $(n$-$1)$-trivial.''
if there exists a projection of $k$, 
such that one can choose $n$ crossings of the projection with the property 
that  making the
crossing changes corresponding to any of the $2^{n}-1$
nontrivial combination of the selected crossings
turns the original knot into the
unknot.  The collection of 
crossing changes is said to be a ``strong $(n$-$1)$-trivializer for $k$''.
\end{definition}
{\em Note:} \hspace{.5mm} The expression ``$n$ adjacent 
to the unknot'' is used interchangeably with ``strongly $(n$-$1)$-trivial.''
We will stick with the latter throughout this
paper.

\medskip
In Section~\ref{sect:bigns} we show that for any $n$ there is a 
non-trivial knot that is strongly $n$-trivial.  On the other hand in 
Section~\ref{sect:onearc} we prove the main result of this paper:

\begin{theorem}
Any non-trivial knot $k$ of  genus $g$
fails to be strongly $n$-trivial for all $n$,
$n \geq 3g-1$.
\label{thm:onearc}
\end{theorem}
{\em Note:} \hspace{.5mm} A knot that is strongly $n$-trivial is 
automatically strongly $m$-trivial for all $m \leq n$.  Also any knot
that is strongly $n$-trivial is clearly $n$-trivial, too.

In analogy with Conjecture~\ref{conj:setsarcs} we have

\begin{corollary}
The unknot is the only knot that is strongly $n$-trivial for all $n$. 
\label{cor:onearc}
\end{corollary}



Theorem~\ref{thm:onearc}
is proven by repeated use of the following theorem of Gabai

\begin{theorem} (Corollary 2.4 [G]) Let $M$ be a Haken manifold
whose boundary is a nonempty union of tori.  Let $F$ be a Thurston
norm minimizing surface representing an element of $H_2(M, \partial M)$
and let $P$ be a component of $\partial M$ such that
$P \cap F = \emptyset$.  Then with at most one exception (up to 
isotopy) $F$ remains norm minimizing in each manifold $M(\alpha)$
obtained by filling $M$ along an essential simple closed curve $\alpha$
in $P$ In particular $F$ remains incompressible in all but at most one
manifold obtained by filling $P$. 
\label{thm:gabai}
\end{theorem}

\section{Notation}
\label{section:notation}
Let $k$ be a knot that is strongly $(n$-$1)$-trivial.  Let
$p:k \rightarrow R^2$ be a projection with crossings $\{a_1, \dots a_n \}$
demonstrating the strong $(n$-$1)$-triviality.  For each $a_i$
let $c_i$ be the small vertical circle that bounds a disk
$D_i$
that intersects $k$ geometrically
twice, but algebraically 0 times.  We call the $c_i$ linking
circles of $k$ and call $D_i$ a crossing disk after 
[ST]. Let $M$ be the link exterior of $k \cup c_1 \cup \dots c_n$
and $P_i$ be the torus boundary component in $M$ corresponding to 
$c_i$. Either $+1$ or $-1$ filling of $P_i$
will result in the desired crossing change depending on
orientation.  We adopt the convention that each $P_i$ will be oriented
so that $+1$ filling of  $P_i$ corresponds to 
the appropriate crossing change dictated by $a_i$.

\section{Irreducibility}
\begin{lemma}
Let $k$ be a nontrivial knot.  Let $\{c_1, \dots c_n \}$ 
be  linking circles for $k$ that show $k$ is strongly $(n$-$1)$-trivial,
then $M$, the exterior of $k \cup c_1 \dots \cup c_n$, is 
irreducible and  $M$ is therefore Haken.
\label{lemma:irreducible}
\end{lemma}
\begin{proof} Assume $M$ is reducible.  Let $S$ be a sphere
that does not bound a ball on either side.  $S$ cannot be
disjoint from $D_1 \cup \dots D_n$ or else it would bound a
ball on the side that does not contain $k$.  Assume $S$ intersects
$D_1 \cup \dots \cup D_n$ minimally and
transversally.  The intersection will consist
of a union of circles.  Let $r$ be one of these circles that is 
innermost on $S$ (any circle that bounds a disk on $S$ disjoint from 
all the other circles of intersection).  Without loss of 
generality assume $r \subset D_1$.  $r$ cannot
be trivial on $D_1 - (D_1 \cap k)$ since $S \cap D_1$ is minimal.
$r$, however must be trivial in $M$ so must divide $D_1$ into
two pieces, one containing both points of $D_1 \cap k$ and the 
other consisting of an annulus running from $r$ to $c_1$.
This disk on $S$ bounded by $r$
shows that $c_1$
bounds a disk in the exterior of $k$. 
This, however, means that $+1$ surgery on $c_1$ 
leaves $k$ unchanged instead of turning it into an unknot, yielding
the desired contradiction.  
\end{proof}

\section{Minimal genus Seifert surfaces}

This section is dedicated to proving the following theorem.

\begin{theorem} If $k$ has a strong $(n$-$1)$-trivializer
$\{ c_1, \dots c_n \}$ and $F$ is a Seifert surface for $k$
disjoint from 
$\{ c_1, \dots c_n \}$ which is minimal
genus among all such surfaces, 
then $F$ is a 
minimal genus Seifert surface for $k$.
\label{thm:surface}
\end{theorem}

\begin{proof}
Because the $c_i$ have linking number 0 with $k$ we can find a 
Seifert surface for $k$ disjoint from the $c_i$.
Let $F$ be a minimal genus Seifert surface for $k$
in the link complement.  

We supplement the notation introduced in Section~\ref{section:notation}. 
Recall $M$ is
the link exterior of $k \cup c_1 \cup \dots \cup c_n$. 
Let $L$ be the corresponding link of $n+1$ components in
$S^3$.  
$P_i$ is the torus boundary component in $M$ corresponding to 
$c_i$.
Let $M(\alpha)$ refer to the manifold
obtained by filling $M$ along an essential simple closed curve 
of slope $\alpha$
in $P_n$. When $\alpha = 1/m, m \in Z$, $M(\alpha)$ is a link exterior.
Let $L_{\alpha}$ be the corresponding link in $S^3$.  Let $k_{\alpha}$
be the image of $k$ in  $L_{\alpha}$ and $F_{\alpha}$ be the image
of $F$  in  $L_{\alpha}$.

We now prove Theorem~\ref{thm:surface} by induction on $n$.
If $F$ is ever a disk then Theorem~\ref{thm:surface}
is clearly true, so \underline{we will assume that $F$ is not a disk throughout} \underline{the 
proof}.

\bigskip

{\bf The base case:}
Let $k$ be a strongly 
0-trivial knot.  This means that $k$ is unknotting number
1 and there is one linking circle $c_1$ that dictates a crossing change
that unknots $k$.

By Lemma~\ref{lemma:irreducible} if $M$ is 
reducible, then $k$ is the unknot.  As in the proof of  
Lemma~\ref{lemma:irreducible} $c_1$ bounds a disk in the complement of
$k$, so $k \cup c_1$ is the unlink on two components. Therefore, 
$F$ being least genus must be a disk, which is a contradiction,
verifying the claim for $M$ reducible and $n=1$.  
We may assume $M$ is irreducible to
complete the base case. $k_1$ is the unknot.  
Since $F_1$ is not a disk, it is no longer
norm minimizing after the filling.  Thus by Theorem~\ref{thm:gabai}
$F$ is norm minimizing under any other filling of $P_1$.  In particular
$F_{\infty}$ is Thurston  norm minimizing for $L_{\infty}$, 
which is just $k$.
Thus, $F$ is a least genus Seifert surface for $k$.

\bigskip

{\bf The inductive step:}  Now we assume that  if $k$ has a 
strong $(n$-$2)$-trivializer
$\{ c_1, \dots c_{(n-1)} \}$ and $F$ is a Seifert surface for $k$
disjoint from 
$\{ c_1, \dots c_{n-1} \}$, which is minimal
genus among all such surfaces, 
then $F$ is also a 
minimal genus Seifert surface for $k$  and show that the same 
must be true for any strong $(n$-$1)$-trivializer for $k$.

Again by Lemma~\ref{lemma:irreducible} if $M$
is reducible, $k$ must be the unknot.  As in previous arguments, 
the separating
sphere $S$ must intersect at least one $D_i$
in a curve that is essential on $D_i - k$.
Without loss of generality, we may assume that $D_n$
is such a disk. Then $c_n$ bounds a disk in the complement of $k \cup \{c_1, \dots
c_{n-2}\}$.
Since $\{ c_1, \dots c_{n-1} \}$ forms a strong $(n$-$2)$-trivializer
for $k$, the induction assumption implies $k$ bounds 
a disk $\Delta$ disjoint from 
$c_1 \cup \dots c_{n-1}$. Since $c_n$ bounds a disk disjoint from 
$k \cup c_1 \cup \dots c_{n-1}$,
$\Delta$ can clearly be chosen to be
disjoint from $c_n$, too, but this contradicts the assumption 
that $F$ was minimal genus, but not a disk.

We now may finish the proof of Theorem~\ref{thm:surface}
knowing that $M$ is irreducible.   
$k_1$ is an unknot in the link $L_1$.
$\{ c_1, \dots c_{n-1} \}$ is a strong $(n$-$2)$-trivializer
for $k$ in $L_1$. The
inductive assumption means that $k_1$ bounds a disk in 
the exterior of $L_1$.  
This disk is in the same
class as $F_1$ in $H_2(M(1), \partial M(1))$, thereby 
showing that $F_1$ is not Thurston
norm minimizing.  Thus, by  Theorem~\ref{thm:gabai}
$F$ remains 
norm minimizing under any other filling of $P_n$. In particular
$F_{\infty}$ is Thurston  norm minimizing in $L_{\infty}$.
Thus, $F$ is a least genus Seifert surface for $k$
in the complement of $\{ c_1, \dots c_{n-1} \}$.
 $\{ c_1, \dots c_{n-1} \}$, however, forms a strong $(n$-$2)$-trivializer
for $k$ in $L_{\infty}$.
By the inductive assumption, $F$ must be Thurston norm 
minimizing for $k$ in the knot complement as well
as the link complement.
\end{proof}

\section{Arcs on a Seifert surface}
We now prove Theorem~\ref{thm:onearc}:
Any non-trivial knot $k$ of  genus $g$
fails to be strongly $n$-trivial for all $n$,
$n \geq 3g-1$.

\begin{proof}

Let $k$ be strongly $n$-trivial with $n$-trivializers $\{c_1, \dots c_{n+1} \}$.
Let $F$ be a minimal genus Seifert surface for $k$ disjoint from 
$\{c_1, \dots c_{n+1} \}$ as in Theorem~\ref{thm:surface}. 
 $F$ has genus $g$.

\label{sect:onearc}
Each linking circle $c_i$ bounds a disk $D_i$ that intersects $F$
in an arc running between the two points of $k \cap D_i$
and perhaps some simple closed curves.  Simple closed curves
inessential in $D_i - k$
can be eliminated since $F$ is incompressible.  Any essential curves
$s_j$ must
be parallel to $c_i$ in $D_i - k$.  These curves can be removed
one at a time using the annulus running from $c_i$ to the
outermost $s_j$ to reroute $F$, decreasing the number of intersections.
Thus, if $F$ is assumed to have minimal intersection with 
each of the $D_i$ then it intersects each one in an arc which we shall
call $a_i$ as in Figure~\ref{fig:linkingcircle}. Each $a_i$ is essential in $F$.
Otherwise $c_i$ would bound a disk disjoint from $F$ and the crossing change along
$c_i$ would fail to unknot $k$.

\begin{figure}[htbp] 
\centerline{\hbox{\epsfig{file=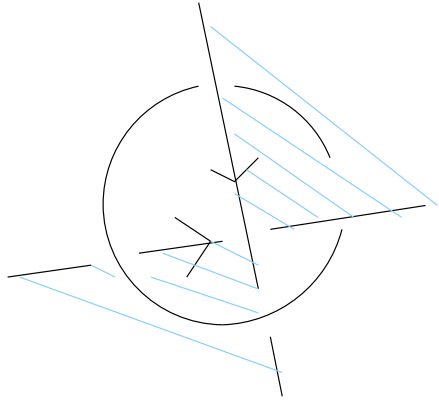}}}
\caption{A Seifert surface passing disjointly through a linking circle}
\label{fig:linkingcircle}
\end{figure}

\begin{lemma}
$a_i$ is never parallel on $F$ to $a_j$ for $i \neq j$.
\end{lemma}
\begin{proof}If $a_i$ is parallel on $F$ to $a_j$
there must be an annulus running from $P_i$ to $P_j$ in the
link exterior.   Recall that we adopted the convention that 
$P_i$ and $P_j$ are each oriented so that $+1$ surgery 
results in the appropriate crossing changes.  
The two tori cannot have opposite
orientations or else $+1$ fillings on both  $P_i$ and $P_j$
is the same as $\infty$ fillings on both, thus,
instead of unknotting $k$  changing
both crossings leaves $k$ knotted.
If the two tori have the same orientation
we could replace $P_i \cup P_j$ with a single
torus $T$ obtained by cutting and pasting
of the two tori along the annulus.
Now $+1$ filling for $P_i$ and $\infty$ filling for $P_j$
is equivalent to $+1$ filling on $T$, while $+1$ filling on 
both $P_i$ and $P_j$ is equivalent to $\frac{1}{2}$ filling on 
$T$. This implies that $F$ fails to be norm-minimizing after
both  $+1$ and  $\frac{1}{2}$
filling of $T$. This contradicts Theorem~\ref{thm:gabai}
completing the proof of the Lemma.
\end{proof}
\vspace{.25in}

Then $\{a_1, \dots a_n \}$ is a collection of 
embedded arcs on $F$, no two of 
which are parallel.  An Euler characteristic argument shows that 
$m \leq 3g-1$.  Since the arcs are in one to one correspondence
with the linking circles,
a strong $n$-trivializer can produce
at most $3g-1$ linking circles for $k$
finally proving Theorem~\ref{thm:onearc}.
\end{proof}

We note that
Theorem~\ref{thm:onearc} predicts that a genus one knot can be at most 
strongly $1$-trivial.
Given standard projections of the trefoil and the figure eight knot
it is easy to find a pair of crossing changes that show the knots are
strongly $1$-trivial.  The theorem is therefore sharp 
at least for genus one knots.  It is possible, but unlikely, that the
theorem remains as precise for higher genus knots.

Finally as noted in the introduction, Theorem~\ref{thm:onearc} implies
Corollary~\ref{cor:onearc}:
The unknot is the only knot that is strongly 
$n$-trivial for all $n$.

\section{Constructing strongly $n$-trivial knots}
\label{sect:bigns}
One might think that
there exists a bound $n$ 
such that no nontrivial knot is strongly $n$-trivial.
Given any $n$, this section gives one way
to construct strongly $n$-trivial knots.

\begin{figure}[htbp] 
\centerline{\hbox{\psfig{figure=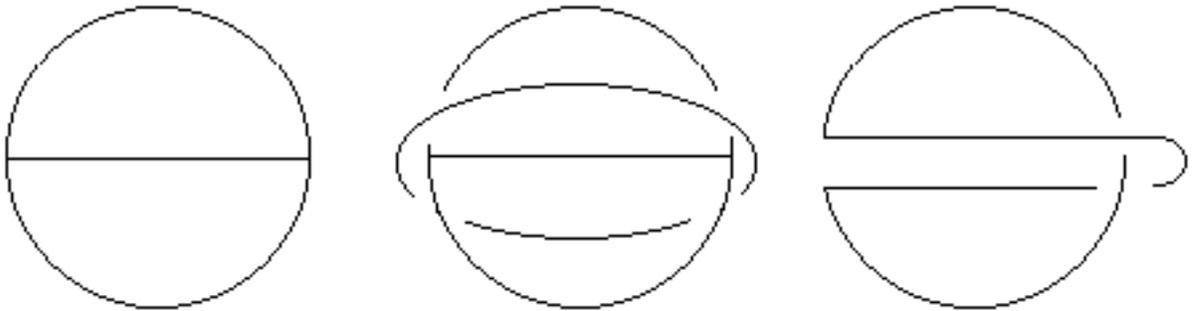}}}
\caption{A graph contains instructions for turning the unknot
into a knot (or perhaps another unknot).}
\label{fig:push}
\end{figure}

\begin{figure}[htbp] 
\centerline{\hbox{\psfig{figure=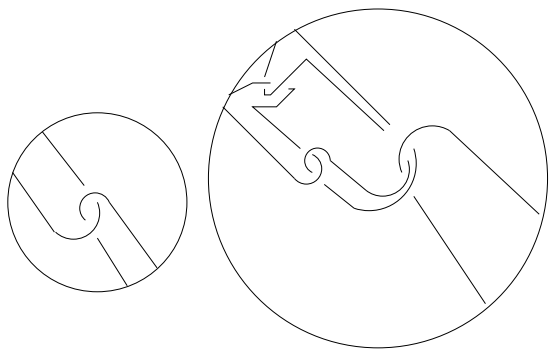}}}
\caption{Examples of crossing changes for the unknot that create
nontrivial knots that are strongly $1$-trivial (left) and 
strongly $2$-trivial (right).  Note that each contains
a subgraph that is a Brunnian link of $n+1$ components.}
\label{fig:sntriv}
\end{figure}

In Figure~\ref{fig:sntriv} we will give 
projections of graphs that show how to turn an unknot into a strongly
$n$-trivial knot.  The circle running around the outside of the 
graph should be viewed as an unknot $k'$.  Each  arc $a_i$ suggests a linking  
circle $c_i$ and a crossing disk $D_i$. If we alter the link $k \cup
c_1 \cup \dots \cup c_n$ in $S^3$ by twisting -1 times
along each of the disks $D_1, D_2, \dots D_n$
$k'$ becomes a new knot $k$
see Figure~\ref{fig:push}.  The linking circles remain fixed,
so we get a new link in $S^3$, $k \cup c_1 \dots \cup c_n$.  

Figure~\ref{fig:sntriv} gives graphs that generate
examples of strongly $1$-trivial 
and strongly $2$-trivial knots.  Note that the figure on the right
is obtained from the figure on the left
by replacing one arc by two new arcs that follow along the original
arc, clasp, return along the original arc, and then, close to the boundary,
clasp once again.  This process could be iterated indefinitely by choosing 
an arc of the new graph and repeating the construction.
It is modeled on doubling one component of a 
link.  Given a Brunnian link of $n$ components 
(a nontrivial link for which any $n-1$ components is the unlink),
doubling one of the components yields a Brunnian link of $n+1$ components.  
The graph on the left in Figure~\ref{fig:sntriv} 
has a Brunnian link of 2 components as a subgraph and 
the one on the right has the double of that link as a subgraph.
Let $\Gamma_{n}$ be the graph after
$n-2$ iterations ($n \geq 2$).

\begin{theorem}
$\Gamma_{n}$ contains a Brunnian link, $l_{n+2}$, of $n+2$ 
components and yields $k$ a non-trivial, strongly $(n+1)$-trivial knot.  
\label{theorem:brunnian}
\end{theorem}

\begin{proof} The link consists of the arcs  $\{a_1, \dots a_{n+2}\}$,
together 
with short segments from $k'$ connecting the end points 
of the segments (and disjoint from the end points of the other segments).  
The base case is trivial because,
$\Gamma_{0}$ contains a Brunnian link of $2$ components: the Hopf link.
$\Gamma_{n}$ is obtained from $\Gamma_{n-1}$ by doubling one of the components
of a Brunnian link of $n+1$ components.  This yields a Brunnian link of $n+2$
components.
 
As a result of the Brunnian structure in $\Gamma_{n}$
any $n+1$ edges from $\{a_1, \dots a_{n+2}\}$ can be 
disjointly embedded on a disk bounded by $k'$. 
So $k'$ forms an unlink with any proper subset of 
$\{c_1, \dots c_{n+2} \}$.
 
We can use this fact to
show that $\{c_1, \dots c_{n+2} \}$ 
are an $n$-trivializer for $k$. 
Let $S$ be any nontrivial subset of $\{c_1, \dots c_{n+2} \}$.
Let $S^c$ be the complement of $S$.
If we take $k$ together with $S$, and do $+1$ surgery on each component
of $S$ the resulting knot is an unknot.  This is because it
is exactly the same as if we took $k'$
and did $-1$ surgery on each of the components of $S^c$.
Since $S$ is a nontrivial subset, $S^c$ is a proper subset.
$k'$ together with a the linking circles in $S^c$, therefore form
an unlink, so each of the components of $S^c$ bounds a disk disjoint
from $k'$ and doing $-1$ surgery on these linking circles
leave $k'$ unchanged.


Now that we know that $k$ is strongly $(n+1)$-trivial,
we need only show $k$ is a non-trivial knot.  
Assume $k$ is trivial. By 
Theorem~\ref{thm:surface}, $k$ bounds a disk $\Delta$
in the complement of $c_1 \cup \dots \cup c_{n+2}$.
Since  $k \cup c_1 \cup \dots
c_{n+2}$ was obtained from $k' \cup c_1 \cup \dots
c_{n+2}$ by spinning along the $D_i$'s, the exteriors of
the two links are homeomorphic, and therefore 
$k'$ must bound a disk $\Delta'$ also disjoint from 
$c_1 \cup \dots \cup c_{n+2}$ (note that one could
even prove that both  $k \cup c_1 \cup \dots \cup c_{n+2}$
and  $k' \cup c_1 \cup \dots \cup c_{n+2}$ are unlinks).
$k'$ intersects
each $D_i$ in 2 points, so as before we may assume $\Delta' \cap D_i$
is an arc for each $i$, but these arcs must, of course,
be isotopic to the $a_is$ which in turn shows that
the $a_is$ can be disjointly embedded on $\Delta'$,
proving that 
$l_{n+2}$ is planar
and not Brunnian, the desired contradiction.
Thus, $k$ is a strongly $(n+1)$-trivial knot, but not the unknot.

\end{proof}

\bigskip

\section{References}

\hspace{.2in}


[G] Gabai, David,  {\it Foliations and the topology of 3-manifolds. II},
Journal of Differential Geometry 26 (1987), pp. 461--478.

\vspace{.21in}

[Gu] Gusarov, M., {\it  On $n$-equivalence of knots and invariants of finite degree}. 
Topology of manifolds and varieties, pp. 173--192, Adv. Soviet Math., 18, Amer. Math. Soc.,
Providence, RI, 1994.

\vspace{.21in}

[NS]  Ng, Ka Yi; Stanford, Ted, {\it On Gusarov's groups of knots}, to appear in Math Proc
Camb Phil. 

\vspace{.21in}

[S] Scharlemann, Martin,  {\it Unknotting number one 
knots are prime}. Invent. Math. 82 (1985), no. 1, 37--55.

\vspace{.21in}


[ST] Scharlemann, Martin; Thompson, Abigail, {\it 
Link genus and the Conway moves}.  Comment. Math.
Helv. 64 (1989), no. 4, 527--535.

\vspace{.21in}

\end{document}